\documentstyle{article}

\newtheorem{theorem}{Theorem}

\newtheorem{definition}[theorem]{Definition}
\newtheorem{example}[theorem]{Example}
\newtheorem{remark}[theorem]{Remark}

\def\QED{\quad\blackslug\lower 8.5pt\null}

\begin{document}

\begin{center}
{\Large \bf  ON CONFORMAL INVARIANCE}
 
\vspace*{2mm}

{\Large \bf  OF ISOTROPIC GEODESICS}

\vspace*{3mm}
{\large M.A. Akivis and  V.V. Goldberg}

\end{center}

\vspace*{5mm}
\begin{abstract}
We consider real isotropic geodesics on manifolds endowed with 
a pseudoconformal structure and their applications to the theory 
of lightlike hypersurfaces on such manifolds, to the geometry of 
four-dimensional conformal structures of Lorentzian type,  and 
to a classification of the Einstein spaces.
\end{abstract}

\setcounter{equation}{0}

\setcounter{section}{0}

{\bf 0. Introduction.}  It is well-known that geodesics on 
a Riemannian or pseudo-Riemannian manifold are not invariant with 
respect to conformal transformations of a Riemannian metric. 
However, in many problems of geometry and physics connected 
with the theory of pseudo-Riemannian manifolds, the isotropic 
geodesics (i.e., geodesics that are tangent to isotropic cones 
at each of their points) arise. 

In the present paper we prove that the isotropic geodesics are 
invariant with respect to conformal transformations of a 
pseudo-Riemannian metric. This is a reason that it is appropriate 
to consider the isotropic geodesics on a manifold endowed 
with a pseudoconformal structure. 

The isotropic geodesics arise naturally when one studies 
lightlike hypersurfaces that are also conformally invariant. 
We prove that such hypersurfaces possess a foliation 
formed by isotropic geodesics. In another terminology 
this fact is given in [DB 96]. We also consider 
the isotropic geodesics on manifolds endowed with 
a pseudoconformal structure of signature (1, 3) and 
apply the obtained results to a classification of the Einstein 
spaces. Note that the geometry of the conformal structures 
of signature $(1, 3)$ was also considered in the recent paper 
[AZ 95].

{\bf 1.} {\bf Preliminaries.} 
Let $V^n_q = (M, g)$ be a pseudo-Riemannian manifold 
of  signature $(p, q)$, where $p+q=n = \dim M$ and $0 < q < n$. 

We  associate with any point 
$x \in M$ its tangent space $T_x (M)$, and define the frame 
bundle whose base is the manifold $M$ and the fibers are the 
families of vectorial frames $\{e_1, \ldots , e_n\}$ in $T_x (M)$ 
defined up to a transformation of the general linear group 
${\bf GL} (n)$.  Let us denote by $\{\omega^1, 
\ldots , \omega^n\}$ the co-frame dual to the frame 
$\{e_1, \ldots , e_n\}$: 
$$
\omega^i (e_j) = \delta^i_j, \;\;\;\;\; i, j = 1, \ldots , n.
$$ 
Then an arbitrary vector $\xi \in T_x (M)$ can be written as 
$$
\xi = \omega^i (\xi) e_i.
$$
The forms $\omega^i$ can be considered as differential forms on 
the manifold $M$ if we assume that $\xi = dx$ is the differential 
of the point $x \in M$. Thus the form $g$ can be 
written as 
\begin{equation}\label{eq:1}
g = g_{ij} \omega^i \omega^j.
\end{equation}
The equation 
\begin{equation}\label{eq:2}
 g_{ij} \omega^i \omega^j = 0
\end{equation}
 defines the foliation of isotropic cones $C_x$ 
of the manifold $M$.

It is well-known that on the manifold $M$ an invariant 
torsion-free affine connection $\gamma$ is defined by means 
of the basis forms $\omega^i$ and connection forms $\omega_j^i$ 
satisfying the following structure equations:
\begin{equation}\label{eq:3}
d \omega^i = \omega^j \wedge \omega_j^i, \;\; d \omega_j^i 
= \omega_j^k \wedge \omega_k^i + R^i_{jkl} \omega^k \wedge 
\omega^l,
\end{equation}
where $R^i_{jkl}$ is the curvature tensor of 
the connection $\gamma$ (see, for example, [C~28]). 
This invariant torsion-free affine connection $\gamma$ 
is called the {\em Levi-Civita connection} if it satisfies 
the equation
\begin{equation}\label{eq:4}
\nabla g_{ij} := dg_{ij} - g_{kj} \omega^k_i - 
g_{ik} \omega_j^k = 0.
\end{equation}
The operator $\nabla$ is called the {\em 
operator of covariant differentiation} with respect to this 
connection. The Levi-Civita connection is uniquely 
determined on $(M, g)$. 

Consider a conformal transformation of a pseudo-Riemannian 
metric, that is, consider another  pseudo-Riemannian metric
\begin{equation}\label{eq:5}
\overline{g} = \sigma g = \sigma g_{ij} \omega^i \omega^j,
\end{equation}
where $\sigma = \sigma (x)$ is a function of a point 
$x \in M$ and $\sigma (x) > 0$. The metric tensor 
of this new metric $\overline{g}$ has the form 
\begin{equation}\label{eq:6}
\overline{g}_{ij} =  \sigma g_{ij}.
\end{equation}
We assume that the basis forms $\omega^i$ are not changed. 

Denote by $\overline{\gamma}$ the Levi-Civita connection 
defined by the metric $\overline{g}$ and by $\overline{\nabla}$ 
the operator of covariant differentiation with respect to the 
connection $\overline{\gamma}$. Then 
\begin{equation}\label{eq:7}
\overline{\nabla} \overline{g}_{ij} =  0.
\end{equation}

It is easy to prove that the connection forms 
$\overline{\omega}_j^i$ of the new 
connection $\overline{\gamma}$ are expressed in terms 
of the connection forms $\omega_j^i$ of the 
connection $\overline{\gamma}$ as follows:
\begin{equation}\label{eq:8}
\overline{\omega}_j^i = \omega_j^i + 
\frac{1}{2} (\delta_j^i \sigma_k \omega^k + 
\sigma_j \omega^i - \sigma^i \omega_j),
\end{equation}
where the quantities $\sigma_k$ are defined by the equation 
\begin{equation}\label{eq:9}
d \log  \sigma = \sigma_k \omega^k
\end{equation}
and 
\begin{equation}\label{eq:10}
\omega_j = g_{jk} \omega^k, \;\; \sigma^j = g^{jk} \sigma_k.
\end{equation}

{\bf 2. Existence  of isotropic geodesics.} Consider geodesics 
on a pseudo-Riemannian manifold $M$. By means of 
 the basis forms $\omega^i$ and connection forms $\omega_j^i$ 
the equations of geodesics on  $M$ can be written in the form
\begin{equation}\label{eq:11}
d \omega^i + \omega^j  \omega_j^i = \alpha \omega^i, 
\end{equation}
where $\alpha$ is a 1-form, and $d$ is the symbol 
of ordinary (not exterior) differentiation.

An {\em isotropic geodesic} on the manifold $M$ 
is a geodesic that is tangent to the isotropic cone $C_x$ 
at each of its points $x$. In addition to equations 
(10), such geodesics also satisfy equation (2) of the isotropic 
cone $C_x$.

We will now prove the following result which confirms 
the existence of isotropic geodesics: 

\begin{theorem}
If a geodesic of the manifold $M$ is tangent to the isotropic 
cone at one of its points $x_0$, then this curve is tangent to 
the isotropic cones at any other of its points $x$; that is, this 
curve is an isotropic  geodesic.
\end{theorem}

{\sf Proof.} The geodesic $x = x (t)$ in question is uniquely 
defined by the system of differential equations (11) 
and initial conditions $x (t_0) = x_0$ and 
$\frac{dx}{dt}|_{t=t_0} = a_0$. Since by hypothesis 
the geodesic is tangent  to the isotropic 
cone at the point $x_0$, we have 
$$
g_{ij}^0 a_0^i a_0^j = 0,
$$
where by $g_{ij}^0$ we denote the values of the components 
of the metric tensor $g_{ij}$ at the point $x_0$. The last 
condition can be rewritten in the form
\begin{equation}\label{eq:12}
(g_{ij} \omega^i \omega^j)|_{x=x_0}  = 0,
\end{equation}
since we have $\omega^i = a^i dt$ along the curve 
$x = x (t)$. 

Differentiating the left-hand side of equation (2) 
and taking into account equations (4) and (11), we 
find that 
$$
d (g_{ij} \omega^i \omega^j) = 2 \alpha g_{ij} \omega^i \omega^j, 
$$
where the differentiation is carried out along the curve 
$x = x (t)$. Along this curve, the 1-form $\alpha$ is a total 
differential and can be written in the form $\alpha 
= d \varphi$, where $\varphi = \varphi (t)$. Thus the last 
equation can be written in the form 
$$
d (g_{ij} \omega^i \omega^j)
= 2 d \varphi \cdot g_{ij} \omega^i \omega^j. 
$$
Integrating this equation, we find that
$$
g_{ij} \omega^i \omega^j = C e^{2\varphi}. 
$$
But since for $t = t_0$ condition (12) holds, we 
find that $C = 0$ and that 
$$
 g_{ij} \omega^i \omega^j = 0
$$
everywhere along the curve $x = x (t)$, so this curve 
is an isotropic geodesic. \rule{3mm}{3mm}

It follows from Theorem 1 that {\em in the pseudo-Riemannian 
manifold $V^n_q, \linebreak q > 0$, through any point $x$ and 
along any isotropic direction emanating from this point, there 
passes one and only one isotropic geodesic}. 

Note that  the usual 
model of space-time in general relativity is a four-dimensional 
pseudo-Riemannian manifold with signature $(1, 3)$ (Lorentzian 
signature) (e.g., see  [Ch~ 83], Ch.~2, \S 11). Isotropic 
geodesics 
of this manifold are curves of propagation of light impulses. 
Hence they are important in this theory.

{\bf 3. Conformal invariance  of isotropic geodesics.} 
On the manifold $V^n_q = (M, \overline{g})$ the equations of 
geodesic lines has the form 
\begin{equation}\label{eq:13}
d \omega^i + \omega^j  \overline{\omega}_i^j = \overline{\alpha} 
\omega^i. 
\end{equation}
Substituting the values (8) of the forms $\overline{\omega}_j^i$ 
into equations (13), we obtain the equations of geodesics 
in the connection $\overline{g}$ in the form
\begin{equation}\label{eq:14}
d \omega^i + \omega^j \omega_i^j - \frac{1}{2} \sigma^i 
g_{jk} \omega^j \omega^k = (\overline{\alpha} 
 - d \log \sigma) \omega^i. 
\end{equation}

Comparing equations (11) and (14), we see 
that {\em in the general case, under conformal  transformation of 
a pseudo-Riemannian metric, geodesics do  not remain invariant}. 
The reason for this is the third term on the left-hand 
side of equation (14) containing $\sigma^i$. However, there 
are two cases where equation (14) defines the same curves 
as equation (11). First of all, this happens if 
$\sigma^i = 0$, that is, if $\sigma = \mbox{const}$. In this 
case equation (14) coincides with  equation (11) with 
$\overline{\alpha} = \alpha$, and all geodesics are transformed 
into geodesics. But this case is not so interesting, since  the 
conformal  transformation has a very special form if 
 $\sigma = \mbox{const}$. Second, equation (14) defines the same 
curves as equation (11) if $g_{jk} \omega^j \omega^k = 0$, 
that is, if the geodesic is isotropic. In this case  
 equations (14) and (11) coincide if 
$\overline{\alpha} = \alpha + d \log \; \sigma$. 
Thus we have proved the following result:

\begin{theorem} 
Under the general  conformal transformation of a 
pseudo-Riemannian metric on a manifold $M$, isotropic geodesics 
and only such geodesics  remain invariant. 
\rule{3mm}{3mm}
\end{theorem}

{\bf 4. The conformal structure \protect\boldmath\(CO (n-\)
\protect\unboldmath 1, 1)} 
The conformal structure 
$CO (n-1, 1)$ on a manifold $M$ of dimension $n$ is a set of 
conformally equivalent pseudo-Riemannian metrics with 
the same signature $(n-1, 1)$. Such a structure is called 
{\em conformally Lorentzian}.

A  metric $g$ can be given on 
$M$ by means of a nondegenerate quadratic form 
$$
g = g_{ij} du^i du^j,
$$
where $u^i, 
i = 1, \ldots , n$, are curvilinear coordinates on  
$M$, and $g_{ij}$ are the components of the metric tensor $g$.

A  pseudoconformal structure on a manifold $M$ is the collection 
of all pseudo-Riemannian metrics  obtained from a fixed  
pseudo-Riemannian metric by conformal transformations. In other 
words, we can say  that a conformal structure on a manifold $M$ is 
defined by means of a relatively invariant quadratic form 
$$
g = g_{ij} du^i du^j.
$$

The equation $g = 0$ defines in the tangent space 
$T_x (M)$ 
a cone $C_x$ of second order called the {\em isotropic 
cone}. Thus the conformal structure $CO (n - 1, 1)$ 
can be given on 
the  manifold $M$ by a field of cones of  second order.

The cone $C_x \subset T_x (M)$ remains invariant under 
transformations of the group 
$$
 G = {\bf SO} (n-1, 1)  \times {\bf H}, 
$$
where ${\bf SO} (n-1, 1)$ is the special $n$-dimensional 
pseudoorthogonal 
group of signature $(n-1, 1)$ (the connected component of the 
unity of the pseudoorthogonal group ${\bf O} (n-1, 1)$), 
and ${\bf H}$ is the group of homotheties. Thus 
the  conformal structure $CO (n-1, 1)$ is a $G$-structure defined 
on the manifold $M$ by the group $G$ indicated above. 
For the conformal structure $CO (n-1, 1)$ the isotropic cone is 
real. 

As in Section {\bf 1}, we  associate with any point 
$x \in M$ its tangent space $T_x (M)$, and define the frame 
bundle whose base is the manifold $M$ and the fibers are the 
families of vectorial frames $\{e_1, \ldots , e_n\}$ in 
$T_x (M)$.  If $\{\omega^1, \ldots , \omega^n\}$ is
 the co-frame dual to the frame $\{e_1, \ldots , e_n\}$, 
then  the form $g$ can be written as 
\begin{equation}\label{eq:15}
g = g_{ij} \omega^i \omega^j
\end{equation}
(see Section {\bf 1}). 

The structure equations of the $CO (n-1, 1)$-structure can be reduced to the following form (see [AG 96], Section {\bf 4.1}):
\begin{equation}\label{eq:16}
d \omega^i = \theta \wedge \omega^i 
+ \omega^j \wedge \theta^i_j, 
\end{equation}
\begin{equation}\label{eq:17}
 d \theta = \omega^i \wedge \theta_i,
\end{equation}
\begin{equation}\label{eq:18}
d \theta^i_j = \theta_j \wedge \omega^i + 
 \theta^k_j \wedge \theta_k^i + g_{jk} \omega^k \wedge 
g^{il} \theta_l + C^i_{jkl} \omega^k \wedge \omega^l,
\end{equation} 
\begin{equation}\label{eq:19}
d \theta_i = \theta_i \wedge \theta + \theta_i^j \wedge \theta_j 
+ C_{ijk}  \omega^j  \wedge  \omega^k. 
\end{equation}
and the metric tensor $g_{ij}$ satisfies the equations
\begin{equation}\label{eq:20}
d g_{ij} - g_{ik} \theta_j^k  -  g_{kj} \theta_i^k  = 0.
\end{equation}
Note that in equations (16)--(20) the forms $\omega^i$ 
are defined in a first-order frame bundle, 
the 1-forms $\theta^i_j$ and a scalar 
 1-form $\theta$ in a second-order frame bundle, and 
a covector form $\theta_i$  in the third-order 
frame bundle. 

For $C^i_{jkl} =  C_{ijk} = 0$,  equations 
(16)--(20) coincide with the structure equations 
 of  the pseudoconformal space $C^n_1$. 
For this  reason  the object $\{C^i_{jkl},  C_{ijk}\}$ is 
called the 
{\em curvature object of the conformal structure $CO (n-1, 1)$}.

 The quantities $C^i_{jkl}$ form a $(1, 3)$-tensor which 
is called the {\em Weyl tensor} or the {\em tensor of conformal 
curvature} of the structure $CO (n-1, 1)$. 

 The quantities $C^i_{jkl}$ and $C_{ijk}$, occurring in equations 
(18) and (19), satisfy the conditions:
\begin{equation}\label{eq:21}
C^i_{jkl} = - C^i_{jlk}, \;\; C_{ijk} = - C_{ikj}, \;\; C^i_{jki} 
= 0.
\end{equation}
These quantities and the tensor $g_{ij}$ also satisfy some 
algebraic and differential equations (see [AG 96], Section 
{\bf 4.1}).

{\bf 5. Isotropic geodesics on lightlike hypersurfaces of 
a manifold \protect\boldmath\(M\;\)\protect\unboldmath endowed 
with \protect\boldmath\(CO (n -\)\protect\unboldmath 
1, 1)-structure.} 
A {\em lightlike hypersurface} $V^{n-1}$ on a manifold 
$M$ of dimension $n$ endowed with a $CO (n-1, 1)$-structure 
of Lorentzian signature $(n-1, 1)$
 is a hypersurface which is tangent to the 
isotropic cone $C_x$ at each point $x \in V^{n-1}$. 

Let $T_x (V^{n-1})$ be the tangent subspace to 
$V^{n-1}$. In $T_x (M)$ we choose a vectorial frame 
$\{e_1, \ldots , e_n\}$ in such a way that its vector $e_1$ has 
the direction of the generator of the isotropic cone $C_x$ along 
which the subspace $T_x (V^{n-1})$ is tangent to $C_x$; 
the vector $e_n$ also has a direction of a generator of the cone 
$C_x$ that does not belong to the subspace $T_x (V^{n-1})$, 
and we locate the vectors $e_a$ in the $(n-2)$-dimensional 
subspace of intersection of $T_x (V^{n-1})$ and the subspace 
tangent to $C_x$ along $e_n$. Then 
\begin{equation}\label{eq:22}
(e_1, e_1) = (e_n, e_n) = 0, \;\; (e_a, e_1) = (e_a, e_n) = 0, 
\;\; (e_a, e_b) = g_{ab}, (e_1, e_n) = - 1,
\end{equation}
where the parentheses denote the scalar product of 
vectors in $T_x (M)$ defined by the quadratic form (15), and 
$a, b = 2, \ldots , n -1$. The last relation in (22) is a 
result of an appropriate normalization of the vectors $e_1$ 
and $e_n$.

With respect to the chosen moving frame,   
 the fundamental form $g$ of $M$ 
can be reduced to the expression
\begin{equation}\label{eq:23}
g = g_{ab} \omega^a \omega^b - 2 \omega^1 \omega^n, 
\;\;\;\;\; a, b = 2, \ldots , n -1, 
\end{equation}
and  the quadratic form $\widetilde{g}$ defining 
the conformal structure on $V^{n-1}$ has the form
\begin{equation}\label{eq:24}
\widetilde{g} = g_{ab} \omega^a \omega^b 
\end{equation}
and is of signature  $(n-2, 0)$, that is, the form 
$\widetilde{g}$ is positive definite. The isotropic cone $C_x$ is 
determined by the equation $g = 0$. 
Thus, the components $g_{ij}$ of the  tensor  $g$ 
are the entries of the following matrix:
\begin{equation}\label{eq:25}
(g_{ij}) = \pmatrix{0 & 0 & -1 \cr 
                                0 & g_{ab} & 0\cr 
                                -1 & 0 & 0}.
\end{equation}

Equations (20) and (25) imply that 
\begin{equation}\label{eq:26} 
\left\{
\renewcommand{\arraystretch}{1.3}
\begin{array}{ll}
\theta_1^n = \theta_n^1 = 0, & \theta_1^1 + \theta_n^n = 0, \\
 \theta_a^n =  g_{ab} \theta^b_1, &
 \theta_a^1 =  g_{ab} \theta^b_n,\\ 
d g_{ab} - g_{ac} \theta_b^c  -  g_{cb} \theta_a^c  = 0. &
\end{array} 
\renewcommand{\arraystretch}{1}
\right.
\end{equation}

Since  the vectors $e_1$ and $e_a$ 
of the frame  $\{e_1, e_a, e_n\}$ of $T_x (M)$ belong to 
$T_x (V^{n-1})$, and the vector $e_n$ 
does not belong to $T_x (V^{n-1})$,  for  
$x \in V^{n-1}$ we have 
$$
dx = \omega_1 e_1 + \omega^a e_a. 
$$
This means that our hypersurface $V^{n-1}$ is defined 
by the following Pfaffian equation:
\begin{equation}\label{eq:27}
\omega_0^n = 0, 
\end{equation}
and the forms $\omega^1$ and $\omega^a, a = 2, \ldots , n - 1,$
 are basis forms of the hypersurface $V^{n-1}$. 

Taking exterior derivative of equation (27) by means of (16), 
we obtain the exterior quadratic equation 
$$
\omega^a \wedge \theta_a^n =0.
$$ 
Applying Cartan's lemma to this equation, 
we find that 
\begin{equation}\label{eq:28}
\theta_a^n = \lambda_{ab} \omega^b, \;\; \lambda_{ab} 
= \lambda_{ba}. 
\end{equation}
It follows from equations (26) and (28) that 
\begin{equation}\label{eq:29}
\theta_1^a = g^{ab} \lambda_{bc} \omega^b,  
\end{equation}
where $g^{ab}$ is the inverse tensor of the tensor 
$g_{ab}$. 

An {\em isotropic geodesic} on the manifold $M$ 
is a geodesic which is tangent to the isotropic cone $C_x$ 
at each of its points $x$. 
We will prove now the following theorem:

\begin{theorem} A  lightlike hypersurface 
$V^{n-1} \subset CO (n-1, 1)$ carries a foliation formed 
 by isotropic geodesics.
\end{theorem}

{\sf Proof}. Since the vectors $e_1$ form an isotropic 
vector field on a hypersurface $V^{n-1}$, 
the equations of the isotropic foliation on $V^{n-1}$ 
have the form
\begin{equation}\label{eq:30}
\omega^a = 0.  
\end{equation}
Thus equation (2) is satisfied. 

Let us prove that the curves belonging to this foliation 
are isotropic geodesics. The isotropic geodesics on the 
conformal structure $CO (n -1, 1)$ are determined by equations 
(2) and (11). However, in Section {\bf 4} we changed 
the notations $\omega_j^i$ for $\theta_j^i$. Thus 
equations (11) can be written as 
\begin{equation}\label{eq:31}
d \omega^i + \omega^j  \theta_j^i = \alpha \omega^i, \;\;\;\;\; 
i, j = 1, \ldots , n,
\end{equation}
where as before $\alpha$ is a 1-form, and $d$ is the symbol of 
ordinary (not exterior) differentiation. 
In our moving frame equations (31) take the form
\begin{equation}\label{eq:32}
\renewcommand{\arraystretch}{1.3}
\begin{array}{ll}
d \omega^1 + \omega^1  \theta_1^1 + \omega^a \theta_a^1 
+ \omega^n  \theta_n^1= \alpha \omega^1, &\\
d \omega^a + \omega^1  \theta_1^a + \omega^b \theta_b^a 
+ \omega^n  \theta_n^a = \alpha \omega^a, &
a, b = 2, \ldots , n-1, \\.
d \omega^n + \omega^1  \theta_1^n + \omega^b \theta_b^n 
+ \omega^n  \theta_n^n = \alpha \omega^n. &
\end{array}
\renewcommand{\arraystretch}{1}
\end{equation}
But by means of equations (27) and (28), which are valid on 
$V^{n-1}$, and equation (30) defining the isotropic 
foliation on $V^{n-1}$, equations (32) are identically 
satisfied. \rule{3mm}{3mm}

Note that a similar result in another terminology is given 
in [DB 96], p. 86.

{\bf 6. Isotropic geodesics on  
a manifold \protect\boldmath\(M\;\;\)\protect\unboldmath  endowed 
with   \protect\boldmath\(CO\)\protect\unboldmath 
 (1, 3)-structure.} 
 We  now consider the pseudoconformal 
$CO (1, 3)$-structure. In a pseudoorthonormal frame its 
fundamental form $g$ becomes
$$
g = - (\omega^1)^2 - (\omega^2)^2 - (\omega^3)^2 
+ (\omega^4)^2.
$$
Transformations of the tangent subspace $T_x (M)$ preserving the  
form $g$ make up the pseudoorthogonal group ${\bf SO} 
(1, 3)$ which is called the {\em Lorentz group}. The isotropic 
cone $C_x \subset T_x (M)$ which is determined by the equation 
$g = 0$ remains invariant under transformations of the group 
$G = {\bf SO} (1, 3) \times {\bf H}$ where ${\bf H}$ is the group 
of homotheties. 

By means of the real transformation
$$
\frac{\omega^1 + \omega^4}{\sqrt{2}} 
\rightarrow \omega^1, \;\; 
\frac{- \omega^1 + \omega^4}{\sqrt{2}} 
\rightarrow \omega^4, \;\; 
\omega^2 \rightarrow \omega^2, \;\; 
\omega^3 \rightarrow \omega^3
$$ 
the form $g$ can be reduced to the form
$$
g = 2 \omega^1 \omega^4 - (\omega^2)^2 - (\omega^3)^2.
$$
It is easy to see now that the cone $g = 0$ carries 
real one-dimensional generators (the straight lines 
$\omega^1 = \omega^2 = \omega^3 
= 0$ and $\omega^4 = \omega^2 = \omega^3 = 0$ are examples of such generators) but does not carry two-dimensional generators.

We complexify   the tangent space $T_x (M)$ by setting 
${\bf C} T_x (M) = T_x (M) \otimes {\bf C}$. 
Moreover, in ${\bf C} T_x (M)$, we will consider only such 
transformations that preserve its real subspace $T_x (M)$ and 
 the symmetry   with respect to $T_x (M)$. 

Next,  by means of the complex transformation 
$$
\omega^1 \rightarrow \omega^1, \;\; 
\frac{\omega^2 +i \omega^3}{\sqrt{2}} 
\rightarrow \omega^2, \;\; 
\frac{\omega^2 - i \omega^3}{\sqrt{2}} 
\rightarrow \omega^3, \;\; 
\omega^4 \rightarrow \omega^4, 
$$
we reduce the quadratic form $g$  to the form 
\begin{equation}\label{eq:33}
g = 2 (\omega^1 \omega^4 - \omega^2 \omega^3),
\end{equation}
where the 1-forms $\omega^1$ and $\omega^4$ are real, and 
the 1-forms $ \omega^2$ and $\omega^3$ are complex conjugate 
forms: 
\begin{equation}\label{eq:34}
 \overline{\omega}^1 =\omega^1, \;\;
 \overline{\omega}^4 = \omega^4, \;\;  
\overline{\omega}^3 = \omega^2.
\end{equation}
It follows that the isotropic cone $g = 0$ carries 
two-dimensional complex conjugate plane generators.

  A vectorial frame in the space 
${\bf C} T_x (M)$, in which  the form $g$ on  the \linebreak 
$CO (1, 3)$-structure  reduces to form (33), 
satisfies the conditions
\begin{equation}\label{eq:35}
\overline{e}_1 = e_1, \;\;  
\overline{e}_4 = e_4, \;\; \overline{e}_2 = e_3.
\end{equation}
Such a frame is called a 
{\em Newman-Penrose tetrad} (see  [NP 62] and 
 [Ch 83], Ch. 1, \S 8). In 
such a frame the vectors $e_1$ and $e_4$ are real, and  the 
vectors $e_2$ and $e_3$ are complex conjugate.

Equations 
\begin{equation}\label{eq:36}
\omega^1 + \lambda \omega^3 = 0, \;\; 
\omega^2 + \lambda \omega^4 = 0
\end{equation}
and 
\begin{equation}\label{eq:37}
\omega^1 + \mu \omega^2 = 0, \;\; 
\omega^3 + \mu \omega^4 = 0
\end{equation}
 determine two families of two-dimensional complex plane 
generators of the complexified cone $C_x$.  
On the $CO (1, 3)$-structure, the parameters 
$\lambda$ and $\mu$ in equations (36) and (37) are complex 
coordinates on the projective lines ${\bf C} P_\alpha$ and 
${\bf C} P_\beta$.  These plane generators 
are, respectively, the {\em $\alpha$-planes} and the 
{\em $\beta$-planes} of the $CO (1, 3)$-structure. If in 
equations (36) we replace all quantities by their conjugates, 
we obtain  equations (37), where $\mu = \overline{\lambda}$. 
Thus there is a one-to-one correspondence between 
$\alpha$-planes and $\beta$-planes  of these two families of 
plane generators of the cone $C_x$, and this correspondence is 
determined by the condition  $\mu = \overline{\lambda}$. 

Since to each point $x \in M$ of a real manifold $M$ carrying a 
 $CO (1, 3)$-{\nolinebreak}struc\-ture 
there correspond two families of 2-planes, 
the family of $\alpha$-planes and the family of $\beta$-planes,  
determined by complex parameters $\lambda$ and $\mu$, two 
bundles, $E_\alpha = (M, {\bf C} P_\alpha)$ and $E_\beta 
= (M, {\bf C} P_\beta)$, arise on $M$, and these two bundles have 
the manifold $M$ as their  common base  and the families of 
complex plane generators of the cone $C_x$ as their fibers. These 
bundles are called the  {\em isotropic bundles} of the 
$CO (1, 3)$-structure.

Since $\mu = \overline{\lambda}$, the isotropic bundles $E_\alpha = (M, {\bf C} P_\alpha)$ and $E_\beta = (M, {\bf C} P_\beta)$ are complex conjugates: 
$\overline{E}_\beta = E_\alpha$.

Two complex conjugate generators of the cone $C_x$ 
determined by the parameters $\lambda$ and 
$\mu = \overline{\lambda}$ 
intersect one another along its real rectilinear generator. 
The equation of this generator can be found from equations 
(36) and (37) provided that $\mu = \overline{\lambda}$. 
Solving these equations, we find that 
$$
\omega^1 = \lambda  \overline{\lambda} \omega^4, \;\; 
\omega^2 = - \lambda \omega^4, \;\; 
\omega^3 =  - \overline{\lambda} \omega^4. 
$$
Hence the directional vector of the  rectilinear generator 
can be written in the form
\begin{equation}\label{eq:38}
\xi =  \lambda  \overline{\lambda} e_1 - \lambda e_2 - 
  \overline{\lambda} e_3 + e_4.
\end{equation} 
Since the basis vectors of the complexified space ${\bf C} T_x$ 
satisfy  relations (35),  
  the vector $\xi$ is real. It depends on 
one complex parameter or two real parameters. Equation (38) 
can be considered as the equation of the director two-dimensional 
surface of the three-dimensional cone $C_x$ in the real 
space $T_x (M)$. 

{\bf 7. Structure equations of the \protect\boldmath\(CO\;\)\protect\unboldmath (1, 3)-structure.} 
In the adapted frame only the following components of the tensor 
$g_{ij}$ will be nonzero: $g_{14} = g_{41} = 1, \; g_{23} 
= g_{32} = - 1$. In view of this,  equations (20) imply that 
the forms $\theta_j^i$ satisfy the conditions 
\begin{equation}\label{eq:39}
\left\{
\renewcommand{\arraystretch}{1.3}
\begin{array}{ll}
\theta_1^4 = \theta_2^3 = \theta_3^2 = \theta_4^1 = 0, \\
\theta_2^4 = \theta_1^3, \;\; \theta_4^2 = \theta_3^1, \;\; \theta_3^4 = \theta_1^2, \;\; \theta_4^3 = \theta_2^1, \\
\theta_1^1 + \theta_4^4 = 0, \;\; \theta_2^2 + \theta_3^3 = 0.
\end{array}
\renewcommand{\arraystretch}{1}
\right.
\end{equation}

\noindent
Now equations (16) take the form
\begin{equation}\label{eq:40}
\left\{
\renewcommand{\arraystretch}{1.3}
\begin{array}{ll}
d\omega^1 = (\theta - \theta_1^1) \wedge \omega^1 + \omega^2 
\wedge \theta_2^1 + \omega^3 \wedge \theta_3^1, \\
d\omega^2 = (\theta - \theta_2^2) \wedge \omega^2 + \omega^1 
\wedge \theta_1^2 + \omega^4 \wedge \theta_3^1, \\
d\omega^3 = (\theta + \theta_2^2) \wedge \omega^3 + \omega^1 
\wedge \theta_1^3 + \omega^4 \wedge \theta_2^1, \\
d\omega^4 = (\theta + \theta_1^1) \wedge \omega^4 + \omega^2 
\wedge \theta_1^3 + \omega^3 \wedge \theta_1^2.
\end{array}
\renewcommand{\arraystretch}{1}
\right.
\end{equation}

By virtue of equations (39), among the forms $\theta_j^i$ 
 only  the forms  $\theta_1^2, \; \theta_2^1, \; \theta_1^3,  
\theta_3^1, \linebreak \theta_1^1$, and $\theta_2^2$ are independent. If 
$\omega^i = 0$, these forms together with the 1-form $\theta$ are 
the invariant forms of a seven-parameter group $G \subset 
{\bf GL} (4)$ that preserves the cone $C_x$ determined by 
equations $g = 0$ where $g$ is determined by equation (33). 

 To write structure equations (18) for 
the $CO (1, 3)$-structure, we consider its tensor of conformal 
curvature $C_{ijkl}$. This tensor has 21 essential nonvanishing 
components that satisfy 11 independent conditions (see [AG 96], 
Section {\bf 5.1}):  
\begin{equation}\label{eq:41}
\left\{
\renewcommand{\arraystretch}{1.3}
\begin{array}{ll}
C_{1234} - C_{1324} + C_{1423} = 0, \\
C_{1224} = C_{1334} = C_{1213} = C_{2434} = 0, \\
C_{1314} - C_{1323} = C_{1424} - C_{2324} = 0, \\
C_{1214} + C_{1223} = C_{1434} + C_{2334} = 0, \\
C_{1414} = C_{2323} = C_{1234} + C_{1324}.
\end{array}
\renewcommand{\arraystretch}{1}
\right.
\end{equation}

Hence the  tensor $C_{ijkl}$  has  10 
independent components in all. We denote them as follows:
\begin{equation}\label{eq:42}
\left\{
\renewcommand{\arraystretch}{1.3}
\begin{array}{lllll}
C_{1212} = a_0, & \!\!\!\! C_{1214} = a_1, & \!\!\!\! 
C_{1234} = a_2,
 & \!\!\!\! C_{1434} = a_3, &\!\!\!\! C_{3434} = a_4, \\
C_{1313} = b_0, &\!\!\!\! C_{1314} = b_1, &\!\!\!\! 
C_{1324} = b_2, 
&\!\!\!\! C_{1424} = b_3, &\!\!\!\! C_{2424} = b_4.
\end{array}
\renewcommand{\arraystretch}{1}
\right.
\end{equation}

\noindent
The remaining components of the  tensor of conformal curvature 
$C_{ijkl}$  are expressible in terms of the above  
components (42).

Now we can write equations (17) and (18) for the 
$CO (1, 3)$-structure in more detail. The former can be written 
as 
\begin{equation}\label{eq:43}
 d \theta = \omega^1 \wedge \theta_1 + \omega^2 \wedge \theta_2 
+ \omega^3 \wedge \theta_3 + \omega^4 \wedge \theta_4,
\end{equation}
and by (39) and (42), the latter has the form 
\begin{equation}\label{eq:44}
\renewcommand{\arraystretch}{1.3}
\begin{array}{ll}
d \theta^1_1 =& \theta_1 \wedge \omega^1 
- \theta_4 \wedge \omega^4  + \theta^2_1 \wedge \theta_2^1 
+ \theta^3_1 \wedge \theta_3^1 \\
&- 2 [a_1 \omega^1 \wedge \omega^2 
+ a_2 (\omega^1 \wedge \omega^4 - \omega^2 \wedge \omega^3) 
+ a_3 \omega^3 \wedge \omega^4 \\
&+ b_1 \omega^1 \wedge \omega^3 
+ b_2 (\omega^1 \wedge \omega^4 + \omega^2 \wedge \omega^3) 
+ b_3 \omega^2 \wedge \omega^4],
\end{array}
\renewcommand{\arraystretch}{1}
\end{equation} 
\begin{equation}\label{eq:45}
\renewcommand{\arraystretch}{1.3}
\begin{array}{ll}
d \theta^2_2 =& \theta_2 \wedge \omega^2 
- \theta_3 \wedge \omega^3  - \theta^2_1 \wedge \theta_2^1 
+ \theta^3_1 \wedge \theta_3^1 \\
&- 2 [a_1 \omega^1 \wedge \omega^2 
+ a_2 (\omega^1 \wedge \omega^4 + \omega^2 \wedge \omega^3) 
+ a_3 \omega^3 \wedge \omega^4 \\
&- b_1 \omega^1 \wedge \omega^3 
- b_2 (\omega^1 \wedge \omega^4 + \omega^2 \wedge \omega^3) 
- b_3 \omega^2 \wedge \omega^4],
\end{array}
\renewcommand{\arraystretch}{1}
\end{equation} 
\begin{equation}\label{eq:46}
\renewcommand{\arraystretch}{1.3}
\begin{array}{ll}
d \theta^2_1 =& \theta_1 \wedge \omega^2 
+ \theta_3 \wedge \omega^4  + (\theta^1_1 - \theta_2^2) 
 \wedge \theta_1^2 \\
&+ 2 [b_0 \omega^1 \wedge \omega^3 
+ b_1 (\omega^1 \wedge \omega^4 + \omega^2 \wedge \omega^3) 
+ b_2 \omega^2 \wedge \omega^4],
\end{array}
\renewcommand{\arraystretch}{1}
\end{equation} 
\begin{equation}\label{eq:47}
\renewcommand{\arraystretch}{1.3}
\begin{array}{ll}
d \theta^1_2 =& \theta_2 \wedge \omega^1 
+ \theta_4 \wedge \omega^3  + \theta^1_2  
 \wedge (\theta_1^1 - \theta_2^2) \\
&- 2 [b_2 \omega^1 \wedge \omega^3 
+ b_3 (\omega^1 \wedge \omega^4 + \omega^2 \wedge \omega^3) 
+ b_4 \omega^2 \wedge \omega^4],
\end{array}
\renewcommand{\arraystretch}{1}
\end{equation} 
\begin{equation}\label{eq:48}
\renewcommand{\arraystretch}{1.3}
\begin{array}{ll}
d \theta^3_1 =& \theta_1 \wedge \omega^3 
+ \theta_2 \wedge \omega^4  
+  (\theta_1^1 + \theta_2^2)  \wedge   \theta^3_1 \\
&+ 2 [a_0 \omega^1 \wedge \omega^2 
+ a_1 (\omega^1 \wedge \omega^4 - \omega^2 \wedge \omega^3) 
+ a_2 \omega^3 \wedge \omega^4],
\end{array}
\renewcommand{\arraystretch}{1}
\end{equation} 

\noindent
and 
\begin{equation}\label{eq:49}
\renewcommand{\arraystretch}{1.3}
\begin{array}{ll}
d \theta^1_3=& \theta_3 \wedge \omega^1 
+ \theta_4 \wedge \omega^2  
+    \theta^1_3  \wedge (\theta_1^1 + \theta_2^2)  \\
&- 2 [a_2 \omega^1 \wedge \omega^2 
+ a_3 (\omega^1 \wedge \omega^4 - \omega^2 \wedge \omega^3) 
+ a_4 \omega^3 \wedge \omega^4].
\end{array}
\renewcommand{\arraystretch}{1}
\end{equation} 

It follows from equations (44) and (45) that 
\begin{equation}\label{eq:50}
\renewcommand{\arraystretch}{1.3}
\begin{array}{ll}
d (\theta^1_1 + \theta^2_2) = &\!\!\!\!2 \theta^3_1 \wedge \theta_3^1 
+ \theta_1 \wedge \omega^1 + \theta_2 \wedge \omega^2  
- \theta_3 \wedge \omega^3 - \theta_4 \wedge \omega^4 \\
&\!\!\!\!- 4 [a_1 \omega^1 \wedge \omega^2 + a_2 (\omega^1 \wedge 
\omega^4 - \omega^2 \wedge \omega^3) + a_3 \omega^3 \wedge 
\omega^4]
\end{array}
\renewcommand{\arraystretch}{1}
\end{equation} 

\noindent
and 
\begin{equation}\label{eq:51}
\renewcommand{\arraystretch}{1.3}
\begin{array}{ll}
d (\theta^1_1 - \theta^2_2) =&\!\!\!\! 
2 \theta^2_1 \wedge \theta_2^1 
+ \theta_1 \wedge \omega^1 - \theta_2 \wedge \omega^2  
+ \theta_3 \wedge \omega^3 - \theta_4 \wedge \omega^4 \\
&\!\!\!\!- 4 [b_1 \omega^1 \wedge \omega^3 + b_2 (\omega^1 \wedge 
\omega^4 + \omega^2 \wedge \omega^3) + b_3 \omega^2 \wedge 
\omega^4].
\end{array}
\renewcommand{\arraystretch}{1}
\end{equation} 

Using   notations (42), we will write now  
10 differential equations  that  the independent 
components of the tensor of conformal curvature $C_{ijkl}$  
satisfy:
\begin{equation}\label{eq:52}
\left\{
\renewcommand{\arraystretch}{1.3}
\begin{array}{ll}
da_0 + 2 a_0 (\theta - \theta_1^1 - \theta_2^2) - 4 a_1
 \theta_1^3 = a_{0i} \omega^i, \\
da_1 +  a_1 (2\theta - \theta_1^1 - \theta_2^2) - a_0 \theta_3^1
-  3 a_2 \theta_1^3 = a_{1i} \omega^i, \\
da_2 + 2 a_2 \theta - 2 a_1 \theta_3^1 - 2a_3 \theta_1^3 
= a_{2i} \omega^i, \\
da_3 +  a_3 (2\theta + \theta_1^1 + \theta_2^2) 
- 3 a_2 \theta_3^1 - a_4 \theta_1^3 = a_{3i} \omega^i, \\
da_4 + 2 a_4 (\theta + \theta_1^1 + \theta_2^2) 
- 4 a_3  \theta_3^1 = a_{4i} \omega^i,
\end{array}
\renewcommand{\arraystretch}{1}
\right.
\end{equation}
\begin{equation}\label{eq:53}
\left\{
\renewcommand{\arraystretch}{1.3}
\begin{array}{ll}
db_0 + 2 b_0 (\theta - \theta_1^1 + \theta_2^2) - 4 b_1
 \theta_1^2 = b_{0i} \omega^i, \\
db_1 +  b_1 (2\theta - \theta_1^1 + \theta_2^2) - b_0 \theta_2^1
- 3 b_2 \theta_1^2 = b_{1i} \omega^i, \\
db_2 + 2 b_2 \theta - 2 b_1 \theta_2^1 - 2b_3 \theta_1^2 
= b_{2i} \omega^i, \\
db_3 +  b_3 (2\theta + \theta_1^1 - \theta_2^2) 
- 3 b_2 \theta_2^1 - b_4 \theta_1^2 = b_{3i} \omega^i, \\
db_4 + 2 b_4 (\theta + \theta_1^1 + \theta_2^2) 
- 4 b_3  \theta_2^1 = b_{4i} \omega^i.
\end{array}
\renewcommand{\arraystretch}{1}
\right.
\end{equation}

\noindent
We can see from  (52) and (53)  that when $\omega^i = 0$, 
the differentials of the components $a_u, \; u = 0, 1, 2, 3, 4$, 
of the tensor of conformal curvature are expressible only in 
terms of these components, and by the same token the same is true 
for the components $b_u$. In view of this,  {\em the tensor of 
conformal curvature of the structure $CO (1, 3)$  is decomposed 
into two subtensors $C_\alpha$ and $C_\beta$ with the components 
$a_u$ and $b_u$, respectively}.

Equations (46)--(51) allow us to establish a geometric 
meaning of the subtensors $C_\alpha$ and $C_\beta$ of the tensor 
of conformal curvature of the $CO (1, 3)$-structure. 
{\em The quantities $a_u$ are the 
components of the curvature tensor of the fiber bundle $E_\alpha$ 
formed by the first family of plane generators of the cones 
$C_x$, 
while the quantities $b_u$ are the components of the curvature 
tensor of the fiber bundle $E_\beta$ formed by the second family 
of plane generators of the cones $C_x$}.

 For the  $CO (1, 3)$-structure not all quantities 
occurring in  equations (39), (40),  and (41)--(49) are real. In 
particular, as we noted earlier, the basis forms $\omega^i$ 
satisfy the equations (34).

The forms  $\theta_j^i$ occurring in equations (40) are 
invariant forms of a complex representation of the real 
six-parameter Lorentz group ${\bf SO} (1, 3)$ that leaves 
invariant the cone $C_x$ determined by the equation $g = 0$ 
in the tangent space $ T_x (M)$. The form  $\theta$ is real,  
$\overline{\theta} =\theta$, since this form is an invariant form 
of the one-parameter group ${\bf H}$ of real homotheties which 
also leaves invariant the cone $C_x$. 

The following theorem is valid (see [AZ 95], Theorem 1): 

\begin{theorem}
On the $CO (1, 3)$-structure,  
the complex forms $\theta_j^i$ occurring in equations $(40)$ 
satisfy the following relations:
\begin{equation}\label{eq:54}
 \overline{\theta}^1_1 = \theta^1_1,  \;\;
 \overline{\theta}^2_2 = - \theta^2_2, \;\;  
\overline{\theta}^3_1 = \theta^2_1, \;\;
 \overline{\theta}^1_3 = \theta^1_2;  
\end{equation}
 the forms $\theta_i$ satisfy the relations
\begin{equation}\label{eq:55}
 \overline{\theta}_1 = \theta_1, \;\;
 \overline{\theta}_2 =  \theta_3, \;\;  
\overline{\theta}_3 = \theta_2, \;\;
 \overline{\theta}_4 = \theta_4;  
\end{equation}
and the components $a_u$ and $b_u, \;\; u = 0, 1, 2, 3, 4$, of 
the curvature tensors $C_\alpha$ and $C_\beta$ of the isotropic 
fiber bundles $E_\alpha$ and $E_\beta$ satisfy the relations 
\begin{equation}\label{eq:56}
\overline{b}_u = a_u.
\end{equation}
\end{theorem}

Let us state some consequences of relations (54)--(56) 
occurring in Theorem 4. 

Equations (54) show that {\em the complex forms $\theta_j^i$ 
occurring in them are expressed in terms of precisely six 
linearly independent real forms}. This number is equal to the 
number 
of parameters on which the Lorentz group depends. These six forms 
are real invariant forms of the group ${\bf SO} (1, 3)$. 

Equations (55) show that among the  forms $\theta_i$ there are 
two real forms and two  complex conjugate forms, and all four 
forms $\theta_i$ are expressed in terms of four linearly 
independent real forms.

Finally, equations (56) show that {\em the  curvature tensors 
$C_\alpha$ and $C_\beta$ of the isotropic fiber bundles  
$E_\alpha$ and $E_\beta$ of the $CO (1, 3)$-structure 
  are  complex conjugates: $\overline{C}_\beta = C_\alpha$}. This 
matches the fact that the  isotropic fiber bundles  
$E_\alpha$ and $E_\beta$ of the $CO (1, 3)$-structure are complex 
conjugates themselves: $\overline{E}_\beta = E_\alpha$. 

It follows that if one of the tensors 
$C_\alpha$ or $C_\beta$ of the $CO (1, 3)$-structure 
vanishes, the other one vanishes too. This implies that {\em the 
$CO (1, 3)$-structure cannot be conformally semiflat 
without being conformally flat}.  
  
{\bf 8. Geometric meaning of the  curvature tensor.} 
 We will now establish a geometric meaning of  the 
subtensors $C_\alpha$ and $C_\beta$ of the curvature tensor of 
the $CO (1, 3)$-structure. Let $\xi = \xi^i e_i$, and let 
$\eta = \eta^i e_i$ be two vectors in the tangent space 
$T_x (M)$, and $\xi \wedge \eta$ be the bivector defined by these 
two vectors. Consider two bilinear forms associated with 
this bivector:
$$
C (\xi \wedge \eta) = C_{ijkl} \xi^i \eta^j \xi^k \eta^l 
=  C_{ijkl} \xi^{[i} \eta^{j]} \xi^{[k} \eta^{l]}
$$
and 
$$
g (\xi \wedge \eta) = (g_{ik} g_{jl} - g_{il} g_{jk}) 
 \xi^i \eta^j \xi^k \eta^l =  (g_{ik} g_{jl} - g_{il} g_{jk}) 
   \xi^{[i} \eta^{j]} \xi^{[k} \eta^{l]}.
$$
Their ratio
$$
K (\xi \wedge \eta) 
= \frac{C (\xi \wedge \eta)}{g (\xi \wedge \eta)} 
$$
is the conformal curvature of the bivector which is called 
the {\em conformal sectional curvature}.

Since $\alpha$-planes and $\beta$-planes are isotropic bivectors, 
for them we have \linebreak 
$g (\xi \wedge \eta) = 0$, and thus the 
expression $K (\xi \wedge \eta)$ does not make sense for them. 
Therefore we will consider for them only the numerator 
$C (\xi \wedge \eta)$ of this expression and will call it the 
{\em relative conformal curvature of two-dimensional isotropic 
direction}.

Let us denote the bivector $\xi \wedge \eta$ by $p$: 
$p = \xi \wedge \eta$, and compute $C (p)$ taking into account 
equations (41) and (42):
\begin{equation}\label{eq:57}
\renewcommand{\arraystretch}{1.3}
\begin{array}{ll}
\displaystyle \frac{1}{4} C (p) 
= & \!\!\!\!  a_0 (p^{12})^2 + 2 a_1 p^{12} (p^{14} 
- p^{23}) + a_2 [2 p^{12} p^{34}  + (p^{14} - p^{23})^2]  \\
&\!\!\!\! + 2 a_3 p^{34} (p^{14} - p^{23}) + a_4 (p^{34})^2 \\
&\!\!\!\! + b_0 (p^{13})^2 + 2 b_1 p^{13} (p^{14} 
+ p^{23}) + b_2 [- 2 p^{13} p^{42}  + (p^{14} + p^{23})^2]  \\
&\!\!\!\! - 2 b_3 p^{42} (p^{14} + p^{23}) + b_4 (p^{42})^2. 
\end{array}
\renewcommand{\arraystretch}{1}
\end{equation}

By (36), the $\alpha$-plane $\alpha (\lambda)$ is determined 
by the vectors
$$
\xi_\lambda = e_3 - \lambda e_1 \;\; \mbox{{\rm and}} \;\;\; 
\eta_\lambda = e_4 - \lambda e_2.
$$
Hence the coordinates of the bivector $p_\lambda = \xi_\lambda 
\wedge \eta_\lambda$ are the minors of the matrix 
$$
\left(
\begin{array}{rrrr}
 - \lambda & 0 & 1 & 0 \\
0 &- \lambda & 0 & 1
\end{array}
\right);
$$
 they are 
$$
p^{12} =  \lambda^2, \;\;  p^{13} = 0, \;\;  
p^{14} = - \lambda, \;\; p^{23} 
= \lambda,\;\;  p^{34} = 1,\;\;  p^{42} = 0.
$$
Substituting these expressions into equations (57), we find 
that 
\begin{equation}\label{eq:58}
\frac{1}{4} C(p_\lambda) = a_0 \lambda^4 
- 4a_1 \lambda^3 + 6 a_2 \lambda^2 - 4 a_3 \lambda + a_4 
:= C_\alpha (\lambda).
\end{equation}

In exactly the same way, by virtue of (37), the $\beta$-plane 
$\beta (\mu)$ is determined by the vectors 
$$
\xi_\mu = e_2 - \mu e_1 \;\; \mbox{{\rm and}} \;\;\; 
\eta_\mu = e_4 - \mu e_3.
$$
This implies that  the coordinates of the bivector 
$p_\mu = \xi_\mu \wedge \eta_\mu$ are 
$$
p^{12} =  0,\;\;  p^{13} = \mu^2, \;\;  p^{14} = - \mu, \;\; 
 p^{23} = - \mu \;\;  p^{34} = 0, \;\; p^{42} = -1,
$$
and the following formula holds:
\begin{equation}\label{eq:59}
 \frac{1}{4}  C(p_\mu) = b_0 \mu^4 
- 4 b_1 \mu^3 + 6 b_2 \mu^2 - 4 b_3 \mu + b_4 := C_\beta (\mu).
\end{equation}

Thus {\em the components of the subtensors $C_\alpha$ and 
$C_\beta$ of the tensor of conformal curvature of a 
$CO (1, 3)$-structure are the coefficients of the polynomials 
$C_\alpha (\lambda)$ and $C_\beta (\mu)$, by means of which we 
can evaluate the relative curvature of the $\alpha$-planes  
$\alpha (\lambda) $ and $\beta$-planes $\beta (\lambda)$, 
respectively.} 

Those isotropic 2-planes  of the structure $CO (1, 3)$ for which  
$C_\alpha (\lambda) = 0$ or $C_\beta (\mu) = 0$ are  called the 
{\em principal $\alpha$-planes} or {\em principal $\beta$-planes} 
of the 
isotropic bundles $E_\alpha$ and $E_\beta$, respectively. Since 
polynomials (58) and (59) are of the fourth degree, it follows 
that, in general, the isotropic cone $C_x$ carries four principal 
$\alpha$-planes and the same quantity of principal  
$\beta$-planes if 
we count each of these planes as many times as its multiplicity. 

By (58) and (59), the equations defining the parameters $\lambda$ 
and $\mu$ of the principal 2-planes of the $CO (1, 3)$-structure 
have the form 
\begin{equation}\label{eq:60}
\left\{
\renewcommand{\arraystretch}{1.3}
\begin{array}{ll}
a_0 \lambda^4 - 4 a_1 \lambda^3 + 6 a_2 \lambda^2 
- 4 a_3 \lambda + a_4 = 0,\\
 b_0 \mu^4 - 4 b_1 \mu^3 + 6 b_2 \mu^2 - 4 b_3 \mu + b_4 = 0.
\end{array}
\renewcommand{\arraystretch}{1}
\right.
\end{equation}
Since by (56)  the 
coefficients of these equations are complex conjugate, 
their solutions $\lambda_p$ and $\mu_p, p = 1, 2, 3, 4,$ 
 are also complex conjugate,  $\mu_p = \overline{\lambda}_p$. 
But as we have proved in Section {\bf 6}, the intersection of 
two complex conjugate 2-planes of the $CO (1, 3)$-structure is a 
real generator of the cone $C_x$. The latter generator is 
determined by the vector 
\begin{equation}\label{eq:61}
\xi_p = \lambda_p  \overline{\lambda}_p e_1 - \lambda_p e_2 
 - \overline{\lambda}_p e_3 + e_4. 
\end{equation}
Thus there arise real fields of principal directions on the 
manifold $M$. 

Now we will prove the following result: 

\begin{theorem}
The integral curves of each of 
four fields of  principal isotropic directions on a manifold $M$ 
endowed with a $CO(1, 3)$-structure are isotropic geodesics of 
the  manifold $M$.
\end{theorem}

{\sf Proof.} As we noted in Sections {\bf 1} and {\bf 2}, 
the isotropic geodesics on a manifold $M$ 
endowed with a $CO(1, 3)$-structure are determined by equations 
(2) and (11). But by (38), in a specialized frame associated with 
a $CO(1, 3)$-structure the coordinates of  
isotropic vectors have the form 
\begin{equation}\label{eq:62}
\xi^1 = \lambda  \overline{\lambda}, \;\; 
\xi^2 = - \lambda, \;\; 
\xi^3 = - \overline{\lambda}, \;\; 
\xi^4 = 1. 
\end{equation}
As a result, equations (11) take the form
\begin{equation}\label{eq:63}
\left\{
\renewcommand{\arraystretch}{1.3}
\begin{array}{ll}
d(\lambda  \overline{\lambda}) - \lambda \theta_2^1 
-  \overline{\lambda} \theta_3^1 =  \lambda  \overline{\lambda} 
(\alpha - \theta_1^1), \\
- d \lambda + \lambda \overline{\lambda}\theta_1^2 
+ \theta_3^1 =  - \lambda  (\alpha - \theta_2^2), \\
- d  \overline{\lambda} +  \lambda \overline{\lambda} \theta_1^3 
+ \theta_2^1 =  -  \overline{\lambda} (\alpha + \theta_2^2), \\
- \lambda \theta_1^3 - \overline{\lambda} \theta_1^2 
 =  \alpha + \theta_1^1.
\end{array}
\renewcommand{\arraystretch}{1}
\right.
\end{equation}

By relations (54),  which the forms $ \theta_j^i$ of 
the $CO (1, 3)$-structure satisfy, only two of equations 
(63), for example, the second and the fourth,  
are independent. Excluding the 1-form $\alpha$ from these  
 equations, we find that 
\begin{equation}\label{eq:64}
d \lambda  + \lambda (\theta_1^1 + \theta_2^2) 
- \theta_3^1 +  \lambda^2 \theta_1^3 = 0.
\end{equation}
Taking exterior derivative of this equation by means 
of (64) and (48)--(50),  we will arrive again to 
equations (60). Moreover, the parameters $\lambda_p$,  
determining the vectors $\xi_p$ in equation (73), 
satisfy equation (64) since for $\lambda = \lambda_p$ 
equations (60) become identities. This means  
that the integral curves of the vector fields $\xi_p$ are  
isotropic geodesics of the  manifold $M$. 
 \rule{3mm}{3mm}

Note also that the integral curves of the principal isotropic 
directions of the $CO (1, 3)$-structure form isotropic geodesic 
congruences on the manifold $M$. In general, the manifold $M$ 
carries four such congruences.
 
{\bf 9. Classification of the Einstein spaces.} 
 For the $CO (1, 3)$-structure, equations (60), 
which by (34) are complex conjugates of one another, 
 are connected with   A.~Z.~Petrov's classification 
 of Einstein spaces. 

We remind 
 that an {\em Einstein space} is a four-dimensional 
pseudo-Riemannian manifold of signature (1, 3)
 whose curvature tensor $R^i_{jkl}$ satisfies the condition 
\begin{equation}\label{eq:65}
R_{jk} - \frac{1}{2} g_{jk} R = - \frac{8 \pi G}{c^4} T_{jk}, 
\end{equation}
where $R_{jk} = R^i_{jki}$ is the Ricci tensor, $R = g^{jk} R_{jk}$ is the scalar curvature of the Riemannian manifold, 
$T_{jk}$ is the energy-momentum tensor, $G$  is the gravitational 
constant, and $c$ is the  speed of light. Equation 
(65) is called the {\em Einstein equation}. 

In empty space, that is, in a region of space-time in which 
$T_{ij} = 0$, the  Einstein equation can be reduced to the form 
$$
R_{ij} = 0.
$$
This implies that the curvature tensor of this space coincides 
with its Weyl 
tensor:   $R^i_{jkl} = C^i_{jkl}$. This follows from the 
expression of the tensor $C^i_{jkl}$ in terms of $R^i_{jkl}, 
R_{jk}$, and $R$ (see [AG 96], Section {\bf 4.2}). 

The classification of  Einstein spaces is connected with the 
structure of its tensor of conformal curvature. Hence this 
classification is of a conformal nature. 
This classification was first constructed by   Petrov 
in   [Pe 54] (see also  [Pi 57]). 

To give a geometric characterization of  Einstein spaces 
of different types, we will also apply the principal 
isotropic congruences 
on the manifolds endowed with a $CO (1, 3)$-structure.

Since for the $CO (1, 3)$-structure, equations (60) 
determining the principal isotropic vector fields $\xi_p$ 
are complex conjugates, for classification of  Einstein spaces 
it is sufficient to consider only one of these equations, for 
example, the first one. By means of this equation, this 
classification can be conducted as follows: 

\begin{enumerate}
\item Type I   of Petrov (we use the Penrose notation for types; 
see [Ch 83], Ch. 1, \S 9, or  [PR 86], Ch.~8) 
is characterized by the fact that all 
roots of equation (60) are distinct. As a result, every isotropic 
cone $C_x$ carries four distinct  principal 
 directions, and the  manifold 
$M$ carries four principal isotropic congruences.

\item Type II of Petrov is characterized by the fact that  
 equation (60) has one double root and two simple roots. 
As a result, every isotropic cone $C_x$ carries three   
principal  directions, one of which is double, 
and the  manifold $M$ carries three 
principal isotropic congruences, one of which is double.

\item Type $D$ of Petrov is characterized by the fact that  
 equation (60) has two distinct double roots. Hence 
every isotropic cone $C_x$ carries two double   
principal  directions, and the  manifold $M$ carries two 
double principal isotropic congruences.

\item Type III of Petrov is characterized by the fact that  
 equation (60) has one triple root and one simple root. 
As a result, every isotropic cone $C_x$ carries two   
principal  directions, one of which is triple,  
and the  manifold $M$ carries two 
 principal isotropic congruences, one of which is triple.

\item Type $N$ of Petrov is characterized by the fact that  
all four roots of equation (60) coincide. 
Hence every isotropic cone $C_x$ carries a quadruple    
principal  direction, and the  manifold $M$ carries a 
quadruple principal isotropic congruence.
\end{enumerate}

{\em Authors' addresses}:\\

\begin{tabular}{ll}
M.A. Akivis &                           V.V. Goldberg \\
Department of Mathematics      &      Department of Mathematics\\
Ben-Gurion University of the Negev & 
                       New Jersey Institute of Technology \\
P.O. Box 653 & University Heights \\
Beer Sheva 84105, Israel & Newark, NJ 07102, U.S.A.     
\end{tabular}

\end{document}